\documentclass{ifacconf}

\usepackage{graphicx}      
\usepackage{natbib}        
\usepackage{amsmath,amssymb,amsfonts}
\usepackage{algorithmic}
\usepackage{graphicx}
\usepackage{textcomp}
\def\BibTeX{{\rm B\kern-.05em{\sc i\kern-.025em b}\kern-.08em
    T\kern-.1667em\lower.7ex\hbox{E}\kern-.125emC}}

\usepackage{nicefrac}       
\usepackage{microtype}      
\usepackage{xcolor}  
\usepackage{mathtools}
 \usepackage{algorithm}
\usepackage{caption, url, balance}
\newcommand{\lmin}{\lambda_{\min}}

\newcommand{\tr}{\textup{trace}}
\newcommand{\bR}{\mathbb{R}}

\newcommand{\MK}{M^{(H)}}
\newcommand{\HK}{J^{(H)}}
\newcommand{\LK}{L^{(H)}}

\newcommand{\citesupp}{}
\newcommand{\citesuppparagraph}{}

\usepackage{xcolor}

\begin{document}
\begin{frontmatter}
\allowdisplaybreaks
\title{On the Relationship of Optimal State Feedback and Disturbance Response Controllers} 

\thanks[footnoteinfo]{The work is
supported by NSF AI institute: 2112085, ONR YIP: N00014-19-1-2217, NSF CNS: 2003111. Contact: runyuzhang@fas.harvard.edu, zhengy@eng.ucsd.edu, weiyuli@g.harvard.edu, nali@seas.harvard.edu}

\author{Runyu (Cathy) Zhang*, Yang Zheng**, Weiyu Li*, Na Li*} 

\address{*School of Engineering and Applied Science, Harvard University}
\address{**Department of Electrical and Computer Engineering, University of California San Diego}

\begin{abstract}                
This paper studies the relationship between state feedback policies and disturbance response policies for the standard Linear Quadratic Regulator (LQR). 
For open-loop stable plants, we establish a simple relationship between the optimal state feedback controller $u_t=K_\star x_t$ and the optimal disturbance response controller $u_t=\LK_{\star;1}w_{t-1}+\ldots+L^{(H)}_{\star;H}w_{t-H}$ with $H$-order. Here $x_t, w_t, u_t$ stands for the state, disturbance, control action of the system, respectively. Our result shows that $L_{\star,1}^{(H)}$ is a good approximation of $K_\star$ and the approximation error $\|K_\star - L_{\star,1}^{(H)}\|$ decays exponentially with $H$. We further extend this result to LQR for open-loop unstable systems, when a pre-stabilizing controller $K_0$ is available.  
\end{abstract}


\end{frontmatter}

\vspace{1mm}

\section{Introduction}

\vspace{-2mm}

Linear quadratic regulator (LQR) is one of the most~fundamental optimal control problems \citep{anderson2007optimal}. Its analytic solution and numerical methods have been well-established in the literature. Specifically, the infinite-time horizon LQR in the discrete time domain is formulated as follows:
\begin{align}
    \!\!\min_{\{u_t\}_{t=0}^{+\infty}} \quad & C\!:=\! \lim_{T\to\infty}\! \frac{1}{T}\mathbb{E}\!\sum_{t=0}^{T-1}\!\!\left(x_t^\top Q x_t \!+\! u_t^\top Ru_t \!+\! 2u_t^\top Sx_t \right)\notag\\
    & \mathrm{subject~to} \quad x_{t+1} = Ax_t + Bu_t + w_t, 
\label{eq:LQR} \vspace{2pt}
\end{align} 
where $x_t \in \mathbb{R}^n$ is the system state, $u_t \in \mathbb{R}^m$ is the control input, $w_t \sim\mathcal{N}(0,I)$ is the external Gaussian process noise, and $Q \succ 0, R \succ 0$ and $S \in \mathbb{R}^{m \times n}$ are performance matrices. Throughout the paper, we make the standard assumption that $
\begin{bmatrix}
    Q & S^\top \\
    S & R
\end{bmatrix}
\succ 0.$ 

It is well-known that the optimal solution for \eqref{eq:LQR} is a state-feedback controller (or policy) $u_t = K_\star x_t$, and the optimal gain $K_\star \in \mathbb{R}^{m \times n}$ can be computed via solving an algebraic Ricatti equation \citep{anderson2007optimal}. The properties of the Ricatti equation and its numerical solutions have been extensively studied \citep{linear-optimal-control-systems, Kleinman1968OnAI,englar1966user}. Most of these results are \textit{model-based }and require the knowledge of system matrices $A, B$ and the weight matrices $Q, R, S$.  
Motivated by the success of model-free policy optimization in reinforcement learning, many recent studies (see review \citet{hu2022towards} and references therein) have started to directly search an optimal policy by viewing the LQR cost $C(K)$ as a function of the policy parameterization $K \in \mathbb{R}^{m \times n}$. This formulation $C(K)$ is more suitable for \textit{model-free policy optimization} but is generally \textit{non-convex}. Thanks to special properties in the optimization landscape such as gradient dominance \citep{fazel2018global}, these methods can still find the optimal controller for the standard LQR problem despite of the non-convexity. However, these properties often fail to generalize to other linear optimal control problems such as sparse or structured LQR, partially observable systems, \citep{zheng2021analysis,hu2022towards}, making it still challenging to develop policy optimization methods with provable convergence and optimality guarantees.

To avoid non-convexity, there are many other methods to re-parameterize the control policy such that the cost function becomes convex under new parameters.~For the general output-feedback case, the classical approach is 
%
%
%
the Youla parameterization \citep{youla1976modern}, and two recent approaches are system-level synthesis (SLS) \citep{wang2019system} and input-output parameterization (IOP) \citep{furieri2019input}; also see \cite{zheng2022system} for two new parameterizations. Another specific idea is to parameterize the control policy as a function of the past disturbances $w_t$, known~as~the \textit{disturbance response control} (DRC) \citep{goulart2006optimization,agarwal2019online}.
In particular, for open-loop stable plants, given a horizon $H \in \mathbb{N}$, we can use a DRC of the form
\begin{equation} \label{eq:DRC}
    u_t=\LK_1w_{t-1}+\ldots+\LK_Hw_{t-H},
\end{equation}
where $\LK_k \in \mathbb{R}^{m \times n}, k = 1, \ldots , H$ are policy parameters, and view the LQR cost in \eqref{eq:LQR} as a function $C(\LK)$ over $\LK:=\{\LK_1, \ldots, \LK_H\}$. It is not difficult to see that the closed-loop state and input evaluations in \eqref{eq:LQR} become affine in $\LK$, and the LQR cost $C(\LK)$ is thus convex in terms of $\LK$. 
%
%
Thanks to the convexity, disturbance-based policy parameterizations 
appear to be easier and more suitable for model-free and online learning setups, which have indeed received increasing attention in online learning and control communities; see e.g., \citet{simchowitz2020improper,li2021safe,agarwal2019online,agarwal2019logarithmic,goulart2006optimization}. It is known that DRC-type controllers are closely related with other convex parameterizations such as the aforementioned Youla, SLS, and IOP.\footnote{For example, interested readers can find some explicit connections in the note: {\scriptsize \url{https://zhengy09.github.io/course/notes/L3.pdf}}.}  

\textbf{Our contribution.} In this paper, we study the relationship 
 between the optimal state feedback policy $K_\star$ and the optimal DRC policy $\LK_\star$ in \eqref{eq:DRC}. 
For open-loop stable plants (i.e., $A$ in \eqref{eq:LQR} is stable), it is not surprising that as the horizon $H$ increases, the optimal performance $C(\LK_\star)$ will improve and converge to the optimal LQR performance $C(K_\star)$. Similar analysis has appeared in \citet{agarwal2019online} but in a slightly different online learning setting. Our paper presents an interesting and not obvious relationship: the first element $\LK_{\star;1}$ in $\LK_\star$ is a good approximation of $K_\star$ and the approximation error decays exponentially with increasing $H$ (Theorem \ref{theorem:relationship-K-L}). Our result points out a possibly simple way of converting disturbance feedback controllers to state feedback controllers -- instead of obtaining a state feedback controller using transfer functions\footnote{That is, solving the transfer function from state $x$ to control $u$ when controller $u$ is in the DRC form \eqref{eq:DRC}.}, we can simply extract $\LK_{\star;1}$ which is already a near optimal state feedback control gain. We further generalize the result to the LQR with an unstable open-loop system through considering DRC with a fixed pre-stabilizing controller $K_0$ (Corollary \ref{coro:unstable-case}).  


The proofs of our results are based on two intuitions:~i) the optimal infinite disturbance response $L_\star^{(\infty)}(z)= \sum_{k=1}^{+\infty} L_{\star;k}^{(\infty)}z^{-k}$ induced by the optimal state feedback $K_\star$ has the exact equivalence $L_{\star,1}^{(\infty) }= K_\star$ (see \eqref{eq:L-infty-def}); ii) as the horizon $H\to+\infty$, the optimal $\LK_{\star;1}$ should converge to $L_{\star;1}^{(\infty)} = K_\star$. In particular, proving (ii) is more technically involved, where we first derive a system of linear equations satisfied by $\LK_\star$ (Lemma \ref{lemma:optimal-L-K}) and then show that $L_\star^{(\infty)}$ is an approximate solution of the linear equations (Lemma \ref{lemma:L-infty-equation} and Corollary \ref{cor:ML+H=Xi}).

\textbf{Related Work.} Some previous studies have built certain relationship between the state representation and other convex parameterizations, e.g., \citep{goulart2006optimization,nett1984connection}. The setting that is most similar to our paper is \citet{goulart2006optimization}, where the authors established an equivalence between the affine state feedback controllers and the affine disturbance feedback controllers. However,  \citet{goulart2006optimization} only considered the finite time horizon problem and dynamic state feedback controllers, which is different from our setting in the infinite-time horizon.  The relationship established in \citet{goulart2006optimization} is very different from our results, and the techniques involved in the proofs are different as well.  

\citesuppparagraph

\vspace{1mm}
\section{Preliminaries and Problem Setup}\label{sec:disturbance-response}
\vspace{-2mm}
In this paper, we consider the infinite-time horizon, time-invariant, discrete time LQR problem as defined in \eqref{eq:LQR}. Throughout this paper, we use $\|\cdot\|$ to denote the matrix $\ell_2$ norm, and $\lmin(X)$ to denote the smallest eigenvalue for a symmetric matrix $X$.



\subsection{State-feedback controllers} 
\vspace{-2mm}
When the plant $(A,B)$ is stabilizable, the optimal solution to problem \eqref{eq:LQR} is a linear state feedback controller $u_t = K_\star x_t$ with (c.f. \cite{anderson2007optimal})
\begin{equation}\label{eq:F-expr}
    K_\star = -(R+B^\top P B)^{-1}(B^\top P A+S),
\end{equation}
where the cost-to-go matrix $P$ satisfies the algebraic Ricatti equation
\begin{equation*}
    P \!=\! A^\top P A - (A^\top PB+S^\top)(R + B^\top PB)^{-1}(B^\top PA+S) + Q.
\end{equation*}

Thus, one natural perspective is to parameterize the policy using a single feedback matrix $K \in \mathbb{R}^{m \times n}$, i.e., $u_t = Kx_t$, which we call as the state feedback representation.
%
As stated in the introduction, this state feedback controller is easy to implement, yet it has one drawback that the LQR cost $C(K)$ becomes non-convex with respect to $K$. 


\subsection{Disturbance response controllers}
\vspace{-2mm}
Another approach to solve the LQR problem \eqref{eq:LQR} is from a disturbance response perspective, which converts the problem to a convex optimization. 

\noindent \textbf{Open-loop stable systems.}
For open-loop asymptotically stable systems, i.e., the spectral radius of $A$ is smaller than $1$,  we can consider a disturbance response controller (DRC) of the form \citep{agarwal2019online,simchowitz2020improper}: 
\begin{equation*}
    u_t = L_1^{(H)} w_{t-1} + \cdots+L_H^{(H)}w_{t-H},
\end{equation*}
where $w_s = 0$ for $s<0$. We can view the LQR cost as a function $C(\LK)$ of the DRC matrices 
$$
\LK:=\{\LK_1, \ldots, \LK_H\}.
$$
We then solve the following optimization problem to get the optimal DRC controller:
\begin{align}
    \min_{\LK} \;\; C(\LK) \!&:= \!\lim_{T\to\infty}  \nonumber \frac{1}{T}\mathbb{E}\sum_{t=0}^{T-1}x_t^\top Q x_t \!+ \! u_t^\top Ru_t \! + \! 2u_t^\top Sx_t  \nonumber\\
    \mathrm{subject~to}\;\; x_{t+1} &= Ax_t + Bu_t + w_t, \label{eq:LQR-disturbance-feedback} \\ 
    u_t &= L_1^{(H)} w_{t-1} + \cdots+L_H^{(H)}w_{t-H}.  \nonumber 
\end{align}
It is not difficult to see that \eqref{eq:LQR-disturbance-feedback} is a convex problem over $\LK$ since the closed-loop state $x_t$ and input $u_t$ all become affine in $\LK$.  






In this paper, we are interested in establishing the relationship between  between the optimal state feedback policy $K_\star$ from \eqref{eq:F-expr} and the optimal DRC policy $\LK_\star$ from \eqref{eq:LQR-disturbance-feedback}. 
First of all, it is not surprising that as the horizon $H$ increases, the optimal performance $C(\LK_\star)$ will improve and converge to the optimal LQR performance $C(K_\star)$. Similar analysis has appeared in \citet{agarwal2019online} but in a slightly different online learning setting. For the self-completeness, we provide our own analysis for the  LQR problem on how $C(\LK_\star)$ approximates  $C(K_\star)$ as $H$ increases in Appendix \ref{apdx:DRC-performance-difference} \citesupp. 
In addition to this relationship between $\LK_\star$ and $K_\star$, we will establish an interesting and not obvious relationship which directly connects the first element $\LK_{\star; 1}$ in $\LK_\star$  with $K_\star$, which will be presented in Theorem~\ref{theorem:relationship-K-L}. 

\noindent\textbf{Open-loop unstable systems.} The above open-loop stability assumption is common for DRC type of controllers, e.g., \citep{agarwal2019logarithmic,simchowitz2020improper}. 
Our results can be easily extended to the unstable case by adding a fixed pre-stabilizing state feedback gain $K_0$ to the DRC, as presented below. 
For unstable system, instead of considering a DRC as in \eqref{eq:DRC}, we consider the following modified DRC with a fixed pre-stabilizing state feedback control gain $K_0$:
\begin{equation}\label{eq:DRC-prestabilizing-K0}
    u_t = K_0x_t+ L_1^{(H)} w_{t-1} + \cdots+L_H^{(H)}w_{t-H}.
\end{equation}
Note that $K_0$ in \eqref{eq:DRC-prestabilizing-K0} is a pre-fixed matrix and does not change when optimizing $C(\LK)$. Given that $K_0$ stabilizes the system, i.e., $A+BK_0$ is stable, we could re-formulate equation \eqref{eq:LQR} by defining an auxiliary variable 
$$\bar u_t := u_t - K_0x_t,$$
then the LQR problem could be re-formulated as:
\begin{align}
    \min_{\LK} C(\LK)&:=\lim_{T\to\infty} \frac{1}{T}\mathbb{E}\sum_{t=0}^{T-1}x_t^\top \bar Q x_t + \bar u_t^\top  R \bar u_t + 2\bar u_t^\top \bar Sx_t \nonumber\\
   \mathrm{s.t.}~~ x_{t+1} &= \bar Ax_t + B\bar u_t + w_t,\label{eq:LQR-reformulate} 
   \\
   \bar{u}_t &= L_1^{(H)} w_{t-1} + \cdots+L_H^{(H)}w_{t-H}, \nonumber
\end{align}
where
\begin{equation}
\begin{split}\label{eq:bar-defs}
    \bar{A} &:= A\!+\!BK_0,  \quad \bar{S} := RK \!+\!S\\
    \bar{Q} &:= Q \!+\! K_0^\top S \!+\! S^\top K_0 \!+\! K_0^\top R K_0.
\end{split}
\end{equation}
Note that $\bar{A} = A+BK_0$ is now a stable matrix. It can also be shown that $\lmin(R\!-\!\bar S\bar Q^{-1}\bar S^\top) >0$; see Lemma \ref{lemma:lmin-ge-0} in the Appendix \citesupp.  
Furthermore, it is not hard to verify that the optimal $\bar{u}_t$ should satisfy $\bar{u}_t = \bar{K}_\star x_t$, where 
$$\bar{K}_\star = K_\star - K_0.$$ Thus by considering the DRC with a pre-stabilizing $K_0$, we could transform the LQR problem \eqref{eq:LQR} with an unstable $A$ to an LQR problem \eqref{eq:LQR-reformulate} with a stable $\bar{A}$. 

\vspace{1mm}
\section{MAIN RESULTS}
\vspace{-2mm}

In this section, we present our main results on the relationship between state feedback policies and disturbance response policies for LQR. 

To characterize the stability degree, we introduce the following definition of exponential stability. 
\begin{defn}[$(\tau,e^{-\rho})$-stability] For $\tau\ge1, \rho>0$, we call a matrix $A$ $(\tau, e^{-\rho})$-stable if
    $\|A^k\|\le \tau e^{-\rho k}.$
\end{defn}


Note that for any open-loop asymptotically stable system, there exist some $\tau\ge1, \rho>0$ such that both $A$ and $A-BK_\star$ are $(\tau, e^{-\rho})$-stable, i.e.,
\begin{equation}\label{eq:open-closed-loop-stable}
    \textstyle \|A^k\|\le \tau e^{-\rho k}, \;\; \|(A-BK_\star)^k\|\le \tau e^{-\rho k}.
\end{equation}
We will use $\tau, \rho$ later in our main result.

\vspace{-1mm}
\subsection{Open-loop stable systems}
\vspace{-2mm}

Our main result in this paper establishes a simple relationship between the optimal control gain $K_\star$ from the algebraic Ricatti equation \eqref{eq:F-expr} and the optimal $L_\star^{(H)}$ from \eqref{eq:LQR-disturbance-feedback}. In particular, we can prove that $\LK_{\star;1}$ is a good approximation of $K_\star$, which is summarized in the theorem below.

\begin{thm}[Main Result]\label{theorem:relationship-K-L} 
For open loop asymptotically stable systems, let $K_\star$ be the optimal feedback gain in \eqref{eq:F-expr}, and  $\LK_\star$ be the optimal solution of \eqref{eq:LQR-disturbance-feedback}. Then, we have
\vspace{2mm}
\begin{equation*}
    \|K_\star - \LK_{\star;1}\| \le \frac{2\tau^3(\|B\|^2\|K_\star\|\|Q\|+\|B\|\|K_\star\|\| S\|)}{\lmin(R\!-\! S Q^{-1} S^\top)(1-e^{-2\rho})^{5/2}}e^{-H\rho},
\end{equation*}
where $\LK_{\star;1}$ denotes the first element in $\LK_\star$. Here $\tau,\rho$ are given in \eqref{eq:open-closed-loop-stable}.
\end{thm}

Theorem \ref{theorem:relationship-K-L} suggests that as long as $H$ is large enough, $L_{\star;1}^{(H)}$ is a good approximation of $K_\star$ and the approximation error decays exponentially w.r.t $H$. Applying certain perturbation analysis arguments (Lemma B.1 in \cite{tu2019gap} and Lemma C.3 in \cite{krauth2019finite}) suggests that for $H$ that is sufficiently large, the system is also closed-loop stable and achieves near-optimal LQR cost with state feedback gain $L_{\star;1}^{(H)}$. Thus instead of implementing the disturbance feedback as 
$$
u_t =  L_{\star;1}^{(H)} w_{t-1} + \cdots+L_{\star;H}^{(H)}w_{t-H}
$$ 
(which is hard to implement because it needs computation and storage of history disturbances $w_{t-k}$), we could simply design a state feedback with gain $L_{\star;1}^{(H)}$, which is much easier to implement and still guarantees near-optimal performance. However, we would also like to point out that Theorem \ref{theorem:relationship-K-L} heavily relies on the fact that the problem is unconstrained. It would be an interesting future direction to figure out whether a similar relationship still holds for constrained or distributed LQ control settings.
\begin{rem}[\textbf{Discussion on the stability assumption}]\label{remark:stability} We would like to emphasize that Theorem \ref{theorem:relationship-K-L} only holds under the open-loop stability assumption, i.e., the spectral radius of $A$ is smaller than 1. Specifically in the proof, one major lemma (Lemma \ref{lemma:optimal-L-K}) will not hold if $A$ is not stable (see more discussion in Remark \ref{rmk:stability} after the lemma). In fact, without the stability assumption, for $H$ that is not large enough, it can be shown that there is no $H$-order DRC that stabilizes the system (see Lemma \ref{lemma:stability-DRC} in the Appendix \citesupp). 
Theorem \ref{theorem:relationship-K-L} also suggests that the approximation error depends on the stability factors $\tau,\rho$, the more stable the system is, the better the approximation error will be.
\end{rem}
\subsection{Extension to unstable systems}\label{sec:extension-unstable}
\vspace{-1mm}

As discussed in the end of Section~\ref{sec:disturbance-response}, we can transform the LQR problem \eqref{eq:LQR} with an unstable $A$ to an LQR problem \eqref{eq:LQR-reformulate} with a stable $\bar{A}$ by considering the DRC with a pre-stabilizing $K_0$: 
\begin{equation*}
    u_t = K_0x_t+ L_1^{(H)} w_{t-1} + \cdots+L_H^{(H)}w_{t-H}.
\end{equation*}

Therefore, we can easily extend Theorem~\ref{theorem:relationship-K-L} to the unstable systems, as shown below,

\begin{cor}[Extension to the unstable case]\label{coro:unstable-case} 
Let $K_\star$ be the optimal feedback gain in \eqref{eq:F-expr}, and $\bar K_\star := K_\star - K_0$,  If both $\|A+BK_0\|$ and $\|A+BK_\star\|$ are $(\tau, e^{-\rho})$-stable, the optimal solution $\LK_\star$ from \eqref{eq:LQR-reformulate} satisfies
\vspace{1mm}
\begin{equation*}
    \|\bar{K}_\star -  \LK_{\star;1}\| \le \frac{2\tau^3(\|B\|^2\|\bar{K}_\star\|\|\bar Q\|+\|B\|\|\bar{K}_\star\|\|\bar S\|)}{\lmin(R\!-\!\bar S\bar Q^{-1}\bar S^\top)(1-e^{-2\rho})^{5/2}}e^{-H\rho},
\end{equation*}
where $\bar S,\bar{Q}$ are defined as in \eqref{eq:bar-defs}. 
\end{cor}

\vspace{1mm}
\section{Proof Sketches}
\vspace{-2mm}

In this section, we present the proof ideas for our main result in Theorem \ref{theorem:relationship-K-L} by a thorough investigation of problem \eqref{eq:LQR-disturbance-feedback}. We first introduce a result from \cite{supp} which shows that the solution to \eqref{eq:LQR-disturbance-feedback} can be explicitly expressed as the solution to a system of linear equations (Lemma \ref{lemma:optimal-L-K}). We next demonstrate that the disturbance response induced by the optimal control gain $K_\star$ is an approximate solution to the linear equations (Lemma \ref{lemma:L-infty-equation}). Combining these two lemmas leads to the final result in Theorem \ref{theorem:relationship-K-L}. 
\subsection{Explicit solution for problem \eqref{eq:LQR-disturbance-feedback}}
\vspace{-2mm}



It is not difficult to see that problem \eqref{eq:LQR-disturbance-feedback} is an unconstrained quadratic  optimization problem w.r.t. the variables $L_1^{(H)},\dots,L_H^{(H)}$. Thus, it is expected that the optimal solution comes from a system of linear equations. Indeed, \cite{supp} has identified these equations, which are formally stated in the following lemma.

\begin{lem}\label{lemma:optimal-L-K} (\cite{supp}) 
For open-loop asymptotically stable systems, the optimal $\LK$ of problem \eqref{eq:LQR-disturbance-feedback} satisfies a set of linear equations 
\begin{equation}\label{eq:L-tilde-eq}
    \MK\LK + \HK=0,
\end{equation}
where $\MK\in\bR^{Hn_u\times Hn_u}$ and $\HK\in\bR^{Hn_u\times n_x}$ are 
\begin{align}
    \hspace{-10pt}\MK\! \!:=\! \begin{bmatrix}
         M_{11}& M_{12}&\cdots &M_{1H} \\
         M_{21}& M_{22}&\cdots &M_{2H}\\
         \vdots&\vdots & &\\
         M_{H1} & M_{H2}&\cdots &M_{HH}
    \end{bmatrix}
    \HK := \begin{bmatrix}
    J_1\\
    \vdots\\
    J_H
    \end{bmatrix}, \label{eq:def-MK-HK}
\end{align}
with submatrices $ M_{km}\in \bR^{n_u\times n_u}, J_k \in \bR^{n_u\times n_x}$ defined as:
\begin{align}
    M_{km}&:=\begin{cases}
        B^\top GB + R, &  k=m\\
        B^\top GA^{k-m}B + SA^{k-m-1}B, &  k> m\\
        B^\top(A^{m-k})^\top GB + B^\top (A^{m-k-1})^\top S^\top, & k<m 
    \end{cases}, \nonumber \\
  J_k&:= B^\top G A^k + SA^{k-1}. \label{eq:def-M-H}
\end{align}
Here $G\in R^{n_x\times n_x}$ is defined as:
\begin{align}\label{eq:def-G}
    G := \sum_{t=0}^\infty (A^t)^\top Q A^t.
\end{align}
\end{lem}
\begin{rem}\label{rmk:stability}
Note that Lemma \ref{lemma:optimal-L-K}  requires $A$ to be exponentially stable;  otherwise the matrix $G$ in \eqref{eq:def-G} is undefined. Since \eqref{eq:LQR-disturbance-feedback} is a quadratic problem with respect to $\LK$, it can be expected that the proof of Lemma \ref{lemma:optimal-L-K} (see \cite{supp}) can be obtained by purely linear algebraic manipulation that writes out $C(\LK)$ explicitly. In the process, there is one step that uses the Taylor series: 
$$(I - z^{-1}A)^{-1} = \sum_{k=0}^{+\infty}z^{-k}A^k,$$ 
which only holds true if $A$ is exponentially stable.
\end{rem}

\subsection{Relationship to optimal state feedback control gain}\label{sec:translate-back-to-state-feedback}
\vspace{-1mm}

We first consider the following disturbance response controller induced by the optimal state feedback gain $K_\star$, which we denoted as $L_\star^{(\infty)}(z)$. That is, $L_\star^{(\infty)}(z)$ is the transfer function from the disturbance signal $\omega$ to the control $u$ when the controller is $u(t)=K_\star x(t)$. It is straightforward to obtain the formulation of $L_\star^{(\infty)}(z)$,
\begin{align}
        L_\star^{(\infty)}(z)& =z^{-1}K_\star(I - z^{-1}(A+BK_\star))^{-1} \notag \\
        &= \sum_{k=1}^{+\infty} L_{\star;k}^{(\infty)}z^{-k},
\end{align}
where 
\begin{equation}
L_{\star;k}^{(\infty)} := K_\star(A+BK_\star)^{k-1},k\ge 1.   \label{eq:L-infty-def}
\end{equation}

Note that implementing the disturbance response controller with transfer function $L_\star^{(\infty)}(z)$ is equivalent to implementing the state feedback controller with control gain $K_\star$. To study the relationship of $\LK_{\star;1}$ and $K_\star$, it is natural to first study the relationship of $\LK_\star$ and $L_\star^{(\infty)}$. We establish the relationship by showing that $L_\star^{(\infty)}$ solves an `infinite dimension' version of equation \eqref{eq:L-tilde-eq} that is  satisfied by $\LK_\star$:
\begin{lem}\label{lemma:L-infty-equation}
The matrices $L_{\star;k}^{(\infty)}$ defined in \eqref{eq:L-infty-def} satisfy
\begin{equation}\label{eq:L-infty-equation}
    \sum_{m=1}^{+\infty}M_{km}L_{\star;m}^{(\infty)}+J_{k} =0,~~\forall~k\ge 1,
\end{equation}
where $M_{km}, J_k$ are defined in \eqref{eq:def-M-H}. 
\end{lem}

Lemma \ref{lemma:L-infty-equation} is the key enabler of proving Theorem \ref{theorem:relationship-K-L}. For structural clearness, we defer the proof of Lemma \ref{lemma:L-infty-equation} to Appendix \ref{apdx:proof-L-infty-equation} \citesupp. We would like to emphasize that the proof of Lemma \ref{lemma:L-infty-equation} is technically involved and may be of independent interest. 

Here we give an intuitive explanation of this lemma. The key insight is that $L^{(\infty)}_\star$ should satisfy an `infinite dimension' version of equation \eqref{eq:L-tilde-eq} (i.e., $H\to+\infty$), which is exactly \eqref{eq:L-infty-equation}. Since $u_t = K_\star x_t$ is globally optimal among all control policies, it is expected that its induced disturbance response $L_\star^{(\infty)}$ solves the optimization problem \eqref{eq:LQR-disturbance-feedback} for $H\to +\infty$. Thus intuitively, it is expected that if we let the horizon $H$ goes to infinity, the solution $L_\star^{(H)}$ for \eqref{eq:L-tilde-eq} will converge to the optimal $L_\star^{(\infty)}$. This is the reason we expect $L_\star^{(\infty)}$~to satisfy \eqref{eq:L-infty-equation}, which is an `infinite dimension' version of \eqref{eq:L-tilde-eq}. The detailed proof is provided in Appendix \ref{apdx:proof-L-infty-equation} \citesupp.

Lemma \ref{lemma:L-infty-equation} immediately results in the following corollary.

\begin{cor}\label{cor:ML+H=Xi}
Define $L_{\star;1:H}^{(\infty)}$ as
\begin{equation}\label{eq:L_star-infty-truncation}
   L_{\star;1:H}^{(\infty)} = \left[
         L_{\star;1}^{(\infty)}\dots;L_{\star;H}^{(\infty)} \right],
\end{equation}
then we have 
\vspace{-5pt}
\begin{equation*}
    \MK L_{\star;1:H}^{(\infty)} + \HK= \mathcal{E},
\end{equation*}
where for all $1\le k\le H$,
\begin{equation*}
    [\mathcal{E}]_k \!=\!\!\!\! \!\sum_{m=H+1}^{+\infty}\!\!\!\!\!B^\top (A^\top)^{m-k-1}(A^\top GB+S^\top)K_\star(A+BK_\star)^{m-1}.
\end{equation*}
\begin{pf}
From Lemma \ref{lemma:L-infty-equation}, we know that
\begin{align*}
    \sum_{m=1}^{+\infty}M_{km}L_{\star;m}^{(\infty)} + J_{k} &= 0,\\
    \Longrightarrow ~~ \sum_{m=1}^{H}M_{km}L_{\star;m}^{(\infty)} + J_k &= -\sum_{m=H+1}^{+\infty}M_{km}L_{\star;m}^{(\infty)} = [\mathcal{E}]_k,
\end{align*}
which completes the proof.
\end{pf}
\end{cor}

\subsection{Proof of Theorem \ref{theorem:relationship-K-L}}
\begin{pf}[of Theorem \ref{theorem:relationship-K-L}]
From Lemma \ref{lemma:optimal-L-K} and Corollary \ref{cor:ML+H=Xi}, we have
\begin{align*}
    &\MK \LK_\star +\HK = 0,~~\MK L_{\star;1:H}^{(\infty)} + \HK = \mathcal{E}. 
\end{align*}
Subtracting these two equations leads to
$$
   \MK(\LK_\star - L_{\star;1:H}^{(\infty)}) = -\mathcal{E}.
   $$
Then, it is not difficult to see that
   \begin{align*}
    &(\LK_\star - L_{\star;1:H}^{(\infty)})^\top (\LK_\star - L_{\star;1:H}^{(\infty)}) = \mathcal{E}^\top (\MK)^{-2} \mathcal{E}\\
  \preceq   & \;\frac{1}{\lmin(\MK)^2}\sum_{k=1}^H[\mathcal{E}]_k^\top[\mathcal{E}]_k\\
  \preceq  & \;  \frac{1}{\lmin(R\!-\!SQ^{-1}S^\top)^2}\sum_{k=1}^H\|[\mathcal{E}]_k\|^2 I,
\end{align*}
where the last inequality uses the result: $\lmin(\MK)\ge \lmin(R\!-\!SQ^{-1}S^\top)$, which can be found in \cite{supp} (Lemma 9). 

We can upper bound the norm of $[\mathcal{E}]_k$ by
\begin{align*}
    \|[\mathcal{E}]_k\|&\le \sum_{m=H+1}^{+\infty}\|B\|^2\|K_\star\|\|GA^{m-k}\|\|(A+BK_\star)^{m-1}\| \\
    &\qquad~~+\|B\|\|K_\star\|\|S\|\|A^{m-k-1}\|\|(A+BK_\star)^{m-1}\|\\
    &\le\frac{\|B\|^2\|K_\star\|\|Q\|\tau^3}{1-e^{-2\rho}}\sum_{m=H+1}^{+\infty}e^{-(2m-k-1)\rho} \\
    &\quad\quad\quad+ \tau^2\|B\|\|K_\star\|\|S\|\sum_{m=H+1}^{+\infty}e^{-(2m-k-2)\rho}\\
    &= \frac{\|B\|^2\|K_\star\|\|Q\|\tau^3}{(1-e^{-2\rho})^2} e^{\!-(2H\!-\!k\!+\!1)\rho} \\
    &\quad\quad\quad\!+\! \frac{\|B\|\|K_\star\|\|S\|\tau^2}{1-e^{-2\rho}}e^{-(2H\!-\!k)\rho},
\end{align*}
where the second inequality uses the result (\citesupp Appendix D)
$$\|GA^m\| \le \frac{\tau^2\|Q\|e^{-\rho m}}{1-e^{-2\rho}}.$$ Thus
\begin{align*}
    &\sum_{k=1}^H\|[\mathcal{E}]_k\|^2\le2 \left(\frac{\|B\|^2\|K_\star\|\|Q\|\tau^3}{(1-e^{-2\rho})^2}\right)^2 \sum_{k=1}^He^{-2(2H-k+1)\rho} \\
    &\qquad+ 2\left(\frac{\|B\|\|K_\star\|\|S\|\tau^2}{1-e^{-2\rho}}\right)^2\sum_{k=1}^He^{-2(2H-k)\rho}\\
    &\le 2\left(\frac{\|B\|^2\|K_\star\|\|Q\|\tau^3}{(1-e^{-2\rho})^2}\right)^2\frac{1}{1-e^{-4\rho}} e^{-2(H+1)\rho}\\
    &\qquad+ 2\left(\frac{\|B\|\|K_\star\|\|S\|\tau^2}{(1-e^{-2\rho})^{3/2}}\right)^2e^{-2H\rho}\\
    &\le 2\left(\frac{\tau^3(\|B\|^2\|K_\star\|\|Q\|+\|B\|\|K_\star\|\|S\|)}{(1-e^{-2\rho})^{5/2}}\right)^2e^{-2H\rho}.
\end{align*}
Finally, we have
\begin{align*}
    &(\LK_{\star;1} \!-\! K_\star)^\top (\LK_{\star;1} \!-\! K_\star) =  (\LK_{\star;1} \!-\! L_{\star;1}^{(\infty)})^\top (\LK_{\star;1} \!-\!  L_{\star;1}^{(\infty)})\\
    &\preceq (\LK_\star - L_{\star;1:H}^{(\infty)})^\top (\LK_\star - L_{\star;1:H}^{(\infty)})\\
    &\preceq \frac{2}{\lmin(R\!-\!SQ^{-1}S^\top)^2}\\
    &\qquad \times \left(\frac{\tau^3(\|B\|^2\|K_\star\|\|Q\|+\|B\|\|K_\star\|\|S\|)}{(1-e^{-2\rho})^{5/2}}\right)^2e^{-2H\rho}I. 
\end{align*}
This leads to 
$$
 \|\LK_{\star;1}\!\!-\!K_\star\|\!\le\! \frac{2\tau^3(\|B\|^2\|K_\star\|\|Q\|\!+\!\|B\|\|K_\star\|\|S\|)}{\lmin(R\!-\!SQ^{-1}S^\top)(1-e^{-2\rho})^{5/2}}e^{\!-\!H\rho},$$
which completes the proof.
\end{pf}

\vspace{1mm}
\section{Numerical examples}

\vspace{-2mm}
We consider the following randomly generated set of system matrices $A,B,Q,R,S$:
\begin{align*}
    A = \left[
    \begin{array}{ccc}
         -0.584 &   0.351  &  0.398\\
   -0.366  & -0.739  &  0.401\\
    0.512 &   0.187 &  -0.761
    \end{array}
    \right]\!\!,~
    &~B = \left[
    \begin{array}{c}
        -0.1659\\
    1.7690\\
   -0.1603
    \end{array}
    \right]\!\!,\\
    Q \!= \!\left[
    \begin{array}{ccc}
          9.549 &  -2.660  &  6.993\\
   -2.660  &  2.702 &  -1.599\\
    6.993  & -1.599  &  8.282
    \end{array}
    \right]\!\!,~
    R =& 2.593,~
    S \!= \!\left[
    \begin{array}{c}
        0.043\\
    0.206\\
   -1.964
    \end{array}
    \right]\!.
\end{align*}

\begin{figure}[t]
    \centering
    \includegraphics[width = .35\textwidth]{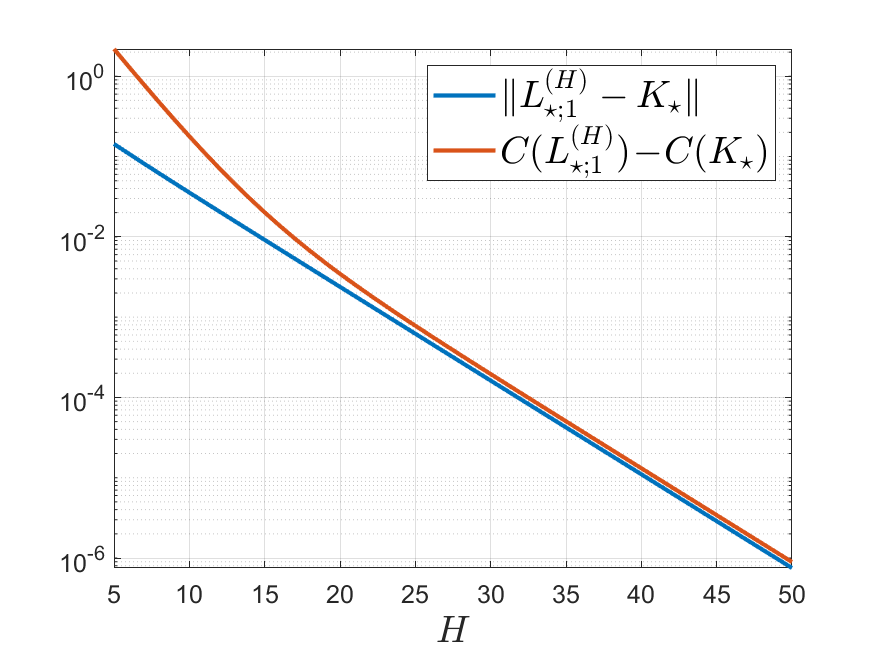}
    \caption{{\small
 Approximation error $\|\LK_{\star;1} - K_\star\|$ and performance difference $C(\LK_{\star;1}) - C(K_\star)$ decays exponentially with $H$}}
    \label{fig:numerics}
\end{figure}

Given the system matrices, we directly call builtin function $\texttt{dlqr}$ in MATLAB System Control Toolbox to solve the optimal state feedback gain $K_\star$. The optimal DRC $\LK_\star$ is solved using \eqref{eq:L-tilde-eq}. Figure \ref{fig:numerics} plots the approximation error $\|\LK_{\star;1} - K_\star\|$ as well as the cost different $C(\LK_{\star;1}) - C(K_\star)$ decays exponentially as $H$ grows larger, which corroborates our theoretical finding in Theorem \ref{theorem:relationship-K-L}.

\vspace{1mm}
\section{Conclusion}
\vspace{-2mm}

This paper has established a simple relationship between the optimal state feedback gain $K_\star$ and the optimal disturbance response controller $u_t= L_{\star;1}^{(H)}w_{t-1}+\ldots+L_{\star; H}^{(H)}w_{t-H}$. The result shows that $L_{\star;1}^{(H)}$ well approximates $K_\star$ and the approximation error decays exponentially with $H$, which points out a possibly simpler way of converting disturbance feedback controllers to state feedback controllers.

\bibliography{bib}

 \appendix
\vspace{-5pt}
\section{Proof of Lemma \ref{lemma:L-infty-equation}}\label{apdx:proof-L-infty-equation}
\vspace{-10pt}
\begin{pf}
Substituting the definitions of $L_{\star;k}^{(\infty)},  M_{km}$ and $J_k$ into the left-hand side of \eqref{eq:L-infty-equation}, we have
\begin{align*}
    &\quad \sum_{m=1}^{+\infty}\!M_{km}L_{\star;m}^{(\infty)} + J_k \\
    &=\! B^\top GA^k \!+\! SA^{k-1} \!+\! \sum_{m=1}^{k-1}\!B^\top GA^{k-m}BK_\star(A+BK_\star)^{m-1} \\
    &+\! \sum_{m=1}^{k-1}\!SA^{k-m-1}BK_\star(A+BK_\star)^{m-1}\\
    &+ \!\sum_{m=k}^{+\infty}B^\top (A^\top)^{m-k} GBK_\star(A+BK_\star)^{m-1}\\
    &+\!\!\!\!\!\sum_{m=k+1}^\infty \!\!\!\!B^\top \!(\!A^\top\!)^{m\!-\!k\!-\!1}S^\top K_\star(A\!+\!BK_\star)^{m\!-\!1}\!\!+\! RK(A\!+\!BK_\star)^{k\!-\!1}\!.
\end{align*}
From the relationship of $K,P$, we have that
\begin{align*}
    &\quad RK_\star(A+BK_\star)^{k-1} \!=\! -(R\!+\!B^\top P B)K_\star(A+BK_\star)^{k-1} \\
    &\qquad\qquad\qquad\qquad\qquad\qquad+ B^\top PBK_\star(A+BK_\star)^{k-1}\\
    &\overset{\eqref{eq:F-expr}}= \!(B^\top P A\!+\!S)(A\!+\!BK_\star)^{k-1} \!\!+\! B^\top PBK_\star(A\!+\!BK_\star)^{k-1} \\
    &= B^\top P(A+BK_\star)^k + S(A+BK_\star)^{k-1},
\end{align*}
which gives
\begin{align*}
    &\quad\sum_{m=1}^{+\infty}M_{km}L_m^{(\infty)} \!+\!J_k \\
    &=\! B^\top \left(GA^k \!+\!\! \sum_{m=1}^{k-1} GA^{k-m}BK_\star(A\!+\!BK_\star)^{m-1}\right)\\
    &\quad +\! S\left(A^{k-1}\! +\!\!\sum_{m=1}^{k-1}A^{k-m-1}BK_\star(A\!+\!BK_\star)^{m-1}\right)\\
     &\quad +B^\top\!\left( \sum_{m=k}^{+\infty} \!(A^\top)^{m-k} GBK_\star(A\!+\!BK_\star)^{m-1}\right.\\
     &\quad+ \!\!\!\left.\sum_{m=k+1}^\infty\!\!(A^\top)^{m\!-\!k\!-\!1}S^\top K_\star(A\!+\!BK_\star)^{m-1} \!+\! P(A\!+\!BK_\star)^k\right) \\
     &\quad -\! S(A\!+\!BK_\star)^{k-1}.
\end{align*}
Since 
\begin{align*}
    &\quad~ GA^k + \sum_{m=1}^{k-1}GA^{k-m}BK_\star(A+BK_\star)^{m-1} \\
    &= GA^k + GA^{k-1}BK_\star + \sum_{m=2}^{k-1}GA^{k-m}BK_\star(A+BK_\star)^{m-1}\\
    &= GA^{k-1}(A+BK_\star) + \sum_{m=2}^{k-1}GA^{k-m}BK_\star(A+BK_\star)^{m-1}\\
    &= G \!\left(\!\!A^{k-1}\!+\!\sum_{m=1}^{k-2}A^{k-1-m}BK_\star(A+BK_\star)^{m-1}\!\right)\!(A\!+\!BK_\star)\\
    &=G\Bigg(A^{k-1}+A^{k-2}BK_\star\\
    &\qquad\qquad +\left.\sum_{m=2}^{k-2}A^{k-1-m}BK_\star(A\!+\!BK_\star)^{m\!-\!1}\right)\!\!(A\!+\!BK_\star)\\
   & =\! G\!\left(\!A^{k-2}\!+\!\sum_{m=1}^{k-3}\!A^{k-1-m}BK_\star(A\!+\!BK_\star)^{m-1}\!\right)\!(A\!+\!BK_\star)^2\\
    &= \cdots\\
    &= GA(A+BK_\star)^{k-1},
\end{align*}
and similarly
\begin{align*}
     S\!\left(\!\!A^{k-1} \!\!+\!\!\sum_{m=1}^{k-1}\!A^{k\!-\!m\!-\!1}BK_\star(A\!+\!BK_\star)^{m\!-\!1}\!\!\right) \!=\! S(A\!+\!BK_\star)^{k-1}
\end{align*}
we can further simplify the expression as
\begin{align*}
    &\quad \sum_{m=1}^{+\infty}M_{km}L_{\star;m}^{(\infty)}+ J_k \\
    &\!=\! B^{\!\top}\!\Big[ \!GA + \!\!\sum_{m=0}^{+\infty} \!(A^{\!\top}\!)^{m} [GBK_\star\!\!\!+\! S^\top K_\star(A\!+\!BK_\star\!)](A\!+\!BK_\star\!)^{m} \\
    &\qquad\qquad\qquad\qquad\qquad\qquad- P(A\!+\!BK_\star) \Big] (A\!+\!BK_\star)^{k-1}.
\end{align*}
\vspace{-5pt}
Let
\begin{equation*}
    X := -\sum_{m=0}^{+\infty} (A^\top)^{m} [GBK_\star+ S^\top K_\star(A+BK_\star)](A+BK_\star)^{m},
\end{equation*}
it suffices to show that $X = GA-P(A+BK_\star)$.
From the definition of $X$ we know that $X$ satisfies the following linear matrix equation
\begin{equation*}
    A^\top X (A+BK_\star) - GBK_\star - S^\top K_\star(A+BK_\star) = X.
\end{equation*}
From the uniqueness of the Sylvester equation (see Lemma \ref{lemma:sylvester}), we know that $X$ is the unique solution to the above linear matrix equation. Thus it suffices to show that
\begin{align*}
     &A^{\!\top} \!(GA\!-\!P(A\!+\!BK_\star\!)) (A\!+\!BK_\star\!) \!-\! GBK_\star \!\!-\! S^{\!\top} K_\star(A\!+\!BK_\star\!)&\\
     &= (GA-P(A+BK_\star))&\\
     &\Longleftrightarrow~~ A^\top GA(A+BK_\star) - A^\top P(A+BK_\star)^2\\
     &\qquad~= G(A+BK_\star)-P(A+BK_\star) + S^\top K_\star (A+BK_\star)\\
     &\Longleftarrow ~~~A^\top GA - A^\top P(A+BK_\star) = G - P + S^\top K_\star.
\end{align*}
From the definition of $G = \sum_{k=0}^{+\infty} (A^\top)^k QA^k$, we have that
\begin{equation*}
    A^\top GA + Q = G,
\end{equation*}
thus it suffices to show that
\begin{align*}
    &- A^\top P(A+BK_\star) = Q - P + S^\top K_\star\\
    \Longleftrightarrow~ &P = A^\top PA + (A^\top PB + S^\top)K_\star + Q\\
    \Longleftrightarrow~ &P = A^\top PA\\
    &~-(A^\top PB + S^\top)(R+B^\top PB)^{-1}(B^\top PA+S) + Q.
\end{align*}
The last equation is exactly the discrete time algebraic Ricatti equation for the optimal cost to go matrix, which completes the proof.
\end{pf}

\vspace{-2pt}
\section{Performance difference}\label{apdx:DRC-performance-difference}
\vspace{-1mm}
In this section we take a deeper look into the relationship of DRC and state feedback control in terms of the performance difference. We show that for any stabilizing state feedback controller $K$, there exists an $H$-order DRC that approximate the cost $C(K)$, where the approximation error decays exponentially with $H$ (Lemma \ref{lemma:DRC-performance-difference}). As a corollary, the performance difference of the optimal DRC and the optimal LQR cost $C(\LK_\star) - C(K_\star)$ also decays exponential with $H$. For any stabilizing $K$, we could define its corresponding equivalent DRC as:
\begin{align*}
    L_K^{(\infty)}(z)& =z^{-1}K(I - z^{-1}(A+BK))^{-1}= \sum_{h=1}^{+\infty} L_{K;h}^{(\infty)}z^{-h},\notag\\
    &\textup{where }  L_{K;h}^{(\infty)} := K(A+BK)^{h-1}, ~~h\ge 1.
\end{align*}
We further define
\begin{equation*}
    L_{K;1:H}^{(\infty)} := [L_{K;1}^{(\infty)};\dots;L_{K;H}^{(\infty)}].
\end{equation*}
We have that the DRC defined by $L_{K;1:H}^{(\infty)}$ has similar cost as $C(K)$, which is formally stated in the following lemma:
\begin{lem}\label{lemma:DRC-performance-difference}
For any $K$ such that $(A+BK)$ is $(\tau,e^{-\rho})$-stable, the $H$-order DRC defined by $L_{K;1:H}^{(\infty)}$ satisfies that
\begin{align*}
   & C(L_{K;1:H}^{(\infty)}) - C(K) \le \\ &n_x^2e^{-2\rho H}\left(\|R\| + \frac{4\tau^4(\|B\|\|K\|^2 \!+\! \|K\|)(\|B\|\|Q\|\!+\!\|S\|)}{(1-e^{-2\rho})^3} \right)
\end{align*}
\end{lem}

Before proving Lemma \ref{lemma:DRC-performance-difference}, we first cite an auxiliary lemma from \cite{supp} that is useful for throughout the proof.
\begin{lem}[\cite{supp}, Appendix B]
\begin{align*}
    C(\LK) = \tr\left(G+ 2{\LK}^\top\HK+{\LK}^\top \MK\LK\right).
\end{align*}
Let $H\to +\infty$ we have that for any $L^{(\infty)}$ that satisfies $\sum_{k}^{+\infty}\|L_k^{(\infty)}\|^2 < +\infty$:
\begin{align*}
    C(L^{(\infty)}) \!=\! \tr\!\left(\!\!G \!+\! 2\sum_{k=1}^{+\infty} {L_k^{(\infty)}}^{\!\top}\! J_k \!+\!\!\!\!\!\!\! \sum_{k=1,m=1}^{+\infty}\!\!\!\! {L_k^{(\infty)}}^{\!\top}\! M_{km} L_m^{(\infty)}\!\!\right)\!.
\end{align*}
\end{lem}
We are now ready to prove Lemma \ref{lemma:DRC-performance-difference}:
\begin{pf}[of Lemma \ref{lemma:DRC-performance-difference}]
\begin{align*}
   &\quad  C(L_{K;1:H}^{(\infty)}) - C(K)=  C(L_{K;1:H}^{(\infty)}) - C(L_{K}^{(\infty)})\\ 
   & = \tr\left(2\sum_{h=H+1}^{+\infty} {L_{K,h}^{(\infty)}}^\top J_h\right.\\&\qquad\left. + \sum_{h=H+1}^{+\infty}\sum_{m=H+1}^{+\infty} {L_{K;h}^{(\infty)}}^\top M_{hm} L_{K;m}^{(\infty)}\right)\\
   &\le n_x^2 \left(2\sum_{h=H+1}^{+\infty} \|L_{K;h}^{(\infty)}\|\| J_h\| \right.\\&\qquad \left.+ \sum_{h=H+1}^{+\infty}\sum_{m=H+1}^{+\infty} \|L_{K;h}^{(\infty)}\| \|M_{hm}\| \|L_{K;m}^{(\infty)}\|\right)\\
   &\le n_x^2 \left(2\sum_{h=H+1}^{+\infty} \|L_{K;h}^{(\infty)}\|\| J_h\| \right.\\&\qquad \left.+ \sum_{h=H+1}^{+\infty}  \left(\sum_{m=H+1}^{+\infty}\|M_{hm}\|\right) \|L_{K;h}^{(\infty)}\|^2
   \right)
\end{align*}
Since $L^{(\infty)}_{K;h} = K(A+BK)^{h-1}, J_h = B^\top GA^h + SA^{h-1}$, from Lemma \ref{lemma:bound-G-norm} we have that
\begin{align*}
    &\|L_{K;h}^{(\infty)}\|\| J_h\| \le \|K\|\left(\frac{\tau e^{-\rho} \|B\|\|Q\|}{1-e^{-2\rho}} + \|S\|\right)\tau^2 e^{-2\rho(k-1)},\\
    &\Longrightarrow~~ \sum_{h=H+1}^\infty\|L_{K;h}^{(\infty)}\|\| J_h\|\\
    &\quad \le \tau^2\|K\|\left(\frac{\tau e^{-\rho} \|B\|\|Q\|}{1-e^{-2\rho}} + \|S\|\right)\sum_{h=H}^{+\infty}e^{-2\rho h}\\
    &\quad = \frac{\tau^2}{1-e^{-2\rho}}\|K\|\left(\frac{\tau e^{-\rho} \|B\|\|Q\|}{1-e^{-2\rho}} + \|S\|\right)e^{-2\rho H}
\end{align*}
Further
\begin{align*}
    &\quad \sum_{m=H+1}^{+\infty}\left\|M_{km}\right\|\\
    &\le \!\|R\|\! +\! \|B^\top GB\| \!+\! 2\!\sum_{m=1}^{+\infty} \!\|B\|^2\|GA^m\| \!+\! 2\!\sum_{m=0}^\infty \!\|B\|\|S\|\|A^m\|\\
    &\overset{\text{(Lemma \ref{lemma:bound-G-norm})}}{\le} \|R\| \!+\! \frac{\tau^2\|B\|^2\|Q\|}{1-e^{-2\rho}} \big( 1\!+\! 2\sum_{m=1}^{+\infty}e^{-\rho m} \big) \!+\! \frac{2\tau\|B\|\|S\|}{1-e^{-\rho}}\\
    &\le \|R\| + \frac{4\tau^2(\|B\|^2\|Q\|+\|B\|\|S\|)}{(1-e^{-2\rho})^2}.
\end{align*}
Moreover,
\begin{align*}
    \sum_{h=H+1}^{+\infty} \|L_{K;h}^{(\infty)}\|^2&\!\le\!\!\!\! \sum_{h=H+1}^{+\infty} \|K\|^2\tau^2 e^{-2\rho (h-1)}\!=\! \frac{\tau^2\|K\|^2}{1-e^{-2\rho}} e^{-2\rho H}\!.
\end{align*}
Combining these bounds together we get:
\begin{align*}
     &\quad C(\LK) - C(K)\le n_x^2e^{-2\rho H} \left(2\frac{\tau e^{-\rho} \|B\|\|Q\|}{1-e^{-2\rho}} + \|S\| \right.\\&+\left.\left(\|R\| + \frac{4\tau^2(\|B\|^2\|Q\|\!+\!\|B\|\|S\|)}{(1-e^{-2\rho})^2} \frac{\tau^2\|K\|^2}{1-e^{-2\rho}}\right) \right) \\
     &\le\! n_x^2e^{-2\rho H}\!\left(\!\|R\|\! +\! \frac{4\tau^4(\|B\|\|K\|^2 \!+\! \|K\|)(\|B\|\|Q\|\!+\!\|S\|)}{(1-e^{-2\rho})^3} \!\right)\!\!.
\end{align*}
\end{pf}

\begin{cor}\label{coro:DRC-performance-difference}
\begin{align*}
    &C(\LK_\star) - C(K_\star) \le \\ &n_x^2e^{-2\rho H}\left(\!\|R\| \!+\! \frac{4\tau^4(\|B\|\|K_\star\|^2 \!+\! \|K_\star\|)(\|B\|\|Q\|\!+\!\|S\|)}{(1-e^{-2\rho})^3} \!\right)\!.
\end{align*}
\begin{pf}
From the optimality of $\LK_\star$, we have that
\begin{align*}
   &C(\!\LK_\star\!) \!-\! C(K_\star) \!=\! C(\!\LK_\star\!) \!-\! C(\!L_\star^{(\infty)}\!) \!\le\!  C(\!L_{\star;1:H}^{(\infty)}\!) \!-\! C(\!L_\star^{(\infty)}\!),
\end{align*}
where $L_{\star;1:H}^{(\infty)} = [L_{\star;1}^{(\infty)}, \dots; L_{\star;H}^{(\infty)}]$. Directly applying Lemma \ref{lemma:DRC-performance-difference} finishes the proof.
\end{pf}
\end{cor}
\vspace{-5pt}
\section{Stability of DRC}
\vspace{-1mm}
\begin{lem}\label{lemma:stability-DRC}
Define $A\in \bR^{n\times n}, B\in \bR^n, e_1\in \bR^n$ as follows:
\begin{align*}
    A = \left[
    \begin{array}{cccc}
         2& 1 && \\
         &\ddots &\ddots&\\
         & &2&1\\
         &&&2
    \end{array}
    \right], B = \left[\begin{array}{c}
         0  \\
         \vdots\\
         0\\
         1
    \end{array}\right], e_1 = \left[\begin{array}{c}
         1  \\
         \vdots\\
         0\\
         0
    \end{array}\right].
\end{align*}
Then for $H\le n$, any DRC of the form \eqref{eq:DRC} is not stable, specifically we have that for $t \ge H$
\begin{align*}
    \mathbb{E} x_{t+1} x_{t+1}^\top \succeq \left(e_1^\top A^H {A^H}^\top e_1\right)\sum_{k=H}^t A^{k-H} e_1 e_1^\top {A^{k-H}}^\top,  
\end{align*}
whose norm blows up exponentially w.r.t. $t$.
\begin{pf}
\begin{align*}
    &\quad x_{t+1} = Ax_t + Bu_t +w_t\\
    &= Ax_t + B(\LK_1w_{t-1}+\ldots+\LK_Hw_{t-H}) + w_t\\
    & = A^2x_{t-1} + w_t + (A + B\LK_1)w_{t-1} + \dots \\
    &\quad~+(AB\LK_{H-1} + B\LK_H) w_{t-H} + AB\LK_H w_{t-H-1}\\
    &=\cdots\\
    &=w_t \!+\! C_1 w_{t-1} \!+\! \dots \!+\! C_{H-1} w_{t-H+1}\! + \!\sum_{k=H}^{t-1}\! A^{k-H} C_H w_{t-k},
\end{align*}
where $C_k$'s are matrices defined by
\begin{align*}
    C_k = A^k + \sum_{t=1}^k A^{k-t}B\LK_t, ~ 1\le k\le H.
\end{align*}
Thus we have that
\begin{align*}
    \mathbb{E} x_{t+1}x_{t+1}^\top &\succeq \sum_{k=H}^{t-1}\! A^{k-H} C_H \left(\mathbb{E}w_{t-k}w_{t-k}^\top\right)C_H^\top{A^{k-H}}^\top\\
    &=\sum_{k=H}^{t-1}\! A^{k-H} C_H C_H^\top{A^{k-H}}^\top.
\end{align*}
Furthermore, since $H \le n$, it is not hard to verify from the definition of $A,B,e_1$ that
\begin{equation*}
    e_1^\top A^{H-k} B = 0, ~1 \le k\le H.
\end{equation*}
Thus
\begin{align*}
    &e_1^\top C_H C_H^\top e_1 = e_1^\top A^H {A^H}^\top e_1\\
    \Longrightarrow~~& C_H^\top C_H \succeq \left( e_1^\top A^H {A^H}^\top e_1\right) e_1e_1^\top.
\end{align*}
Substitute this into the above equation gives
\begin{align*}
    \mathbb{E}x_{t+1}x_{t+1}^\top &\succeq  \left(e_1^\top A^H {A^H}^\top e_1\right)\sum_{k=H}^t A^{k-H} e_1 e_1^\top {A^{k-H}}^\top,
\end{align*}
which completes the proof.
\end{pf}
\end{lem}

\vspace{-5pt}
\section{Auxiliaries}
\vspace{-1mm}

\begin{lem}\label{lemma:bound-G-norm}
For any $m\ge0$, G defined in \eqref{eq:def-G} satisfies
\begin{equation*}
    \|GA^m\| \le \frac{\tau^2\|Q\|e^{-\rho m}}{1-e^{-2\rho}}.
\end{equation*}
\end{lem}
\begin{pf}
From the definition of $G$, we have
\begin{align*}
    &\|GA^m\| = \|\sum_{t=0}^\infty (A^t)^\top Q A^{t+m}\|\le \sum_{t=0}^\infty \|Q\|\|A^t\|\|A^{t+m}\|
    \\&\le \|Q\|\sum_{t=0}^\infty \tau^2 e^{-\rho(2t+m)}=\frac{\tau^2\|Q\|e^{-\rho m}}{1-e^{-2\rho}},
\end{align*}which completes the proof.
\end{pf}

\begin{lem}\label{lemma:lmin-ge-0}
$\lmin(R - \overline{S}\overline{Q}^{-1}\overline{S}^\top) > 0$.
\begin{pf}
It is not hard to check that
\begin{align*}
    &\left[
\begin{array}{ll}
    \overline Q & \overline S^\top \\
    \overline S & R
\end{array}
\right] = \left[
\begin{array}{ll}
    I & K_0^\top \\
    0 & I
\end{array}
\right]\left[
\begin{array}{ll}
    Q & S^\top \\
    S & R
\end{array}
\right]\left[
\begin{array}{ll}
    I & 0 \\
    K_0 & I
\end{array}
\right]\\
&\qquad\Longrightarrow \lmin\left(\left[
\begin{array}{ll}
    \overline Q & \overline S^\top \\
    \overline S & R
\end{array}
\right]\right) \succ 0.
\end{align*}
From Lemma \ref{lemma:schuur-complement-lmin} we have that $\lmin(R - \overline{S}\overline{Q}^{-1}\overline{S}^\top) > 0$.
\end{pf}
\end{lem}

\begin{lem}\label{lemma:schuur-complement-lmin}
$\lmin(R - SQ^{-1}S^\top) \ge \lmin\left(\left[
\begin{array}{ll}
    Q & S^\top \\
    S & R
\end{array}
\right]\right).$
\begin{pf}
Let $\lambda:=\lmin\left(\left[
\begin{array}{ll}
    Q & S^\top \\
    S & R
\end{array}
\right]\right)$, then we have that
\begin{align*}
 &   \left[
\begin{array}{cc}
    Q & S^\top \\
    S & R-\lambda I
\end{array}
\right] \succeq 0 ~~\Longrightarrow~ (R-\lambda I)-SQ^{-1}S^\top \succeq 0\\
&\qquad\qquad \Longrightarrow~R - SQ^{-1}S^\top \succeq \lambda I,
\end{align*}
\end{pf}
which completes the proof.
\end{lem}

\begin{lem}\label{lemma:sylvester}
The discrete time Sylvester equation
\begin{equation}\label{eq:sylvester}
    A^\top X B + C = X, ~~A,B,C\in \bR^{n\times n}
\end{equation}
obtains a unique solution $X\in \bR^{n\times n}$ if and only if any eigenvalue $\lambda_i$ of $A$ and any eigenvalue $\mu_j$ of $B$ satisfies $\lambda_i\mu_j \neq 1$. 
\end{lem} 
\begin{pf}
We denote all the eigenvalues of $A$ as $\lambda_i$'s and all the eigenvalues of $B$ as $\mu_i$'s ($i= 1,2,\dots,n$).
We first prove that if there exists $\lambda_i\mu_j = 1$, then the solution to \eqref{eq:sylvester} is either not unique or has no solution. Since \eqref{eq:sylvester} is a linear equation w.r.t. $X$, it is sufficient to show that
\begin{equation}\label{eq:sylvester-2}
    A^\top X B = X
\end{equation}
has non-zero solutions.
Let $v_i, u_j\in \bR^n$ be nonzero vectors such that
\begin{equation*}
    A^\top v_i = \lambda_i v_i, ~~B^\top u_j = \mu_j u_j,
\end{equation*}
then let $X = v_iu_j^\top$, we have that
\begin{align*}
    A^\top X B = A^\top v_i u_j^\top B = \lambda_i\mu_j v_iu_j^\top = v_iu_j^\top = X,
\end{align*}
which shows that \eqref{eq:sylvester-2} has a non-zero solution.

Next, we prove that if $\lambda_i\mu_j \neq 1$, then \eqref{eq:sylvester} obtains a unique solution. It is sufficient to show that \eqref{eq:sylvester-2}
only have zero solution. From now on we assume that matrix $X$ satisfies $A^\top XB = X$. 
 Denote
\begin{align*}
    p(\lambda) &= \prod_{i=1}^n (\lambda - \lambda_i)\\
    q(\lambda) &= \prod_{i=1}^n (\lambda - \mu_i) = \sum_{k=0}^n b_k \lambda^k\\
    \overline{q}(\lambda) &= \prod_{i=1}^n (1 - \lambda\lambda_i) = \sum_{k=0}^n b_k\lambda^{n-k}.
\end{align*}
From Hamilton-Cayley theorem, we know that
\begin{equation*}
    q(B) = 0,
\end{equation*}
thus
\begin{align*}
    0&=Xq(B) = \sum_{k=0}^n b_k XB^k\\
    &=\sum_{k=0}^n b_k (A^\top)^{n-k}XB^n~~\quad {\left((A^\top)^{n-k}XB^{n-k} = X\right)}\\
    & = \overline{q}(A)XB^n.
\end{align*}
Since $\lambda_i\mu_j \neq 1$, we have that
$\overline{q}(\lambda_i)\neq 0, \forall 1\le i\le n$. Thus, there exists polynomials $f(\lambda), g(\lambda)$ such that
\begin{equation*}
    f(\lambda)p(\lambda) + g(\lambda)\overline q(\lambda) = 1.
\end{equation*}
Thus
\begin{align*}
    f(A)p(A) + g(A)\overline q(A) = I
    \Longrightarrow~ g(A)\overline{q}(A) = I\\
    (p(A) = 0 \textup{ from Hamilton-Cayley theorem}).
\end{align*}
Thus we have
\begin{equation*}
    0 = g(A) \overline{q}(A)XB^n = XB^n.
\end{equation*}
Additionally, we have
\begin{equation*}
    X = (A^\top)^n X B^n = 0,
\end{equation*}
which completes the proof.
\end{pf}

\end{document}